\input amstex
\documentstyle{amsppt}
%
\catcode`@=11
\redefine\output@{%
  \def\break{\penalty-\@M}\let\par\endgraf
  \ifodd\pageno\global\hoffset=105pt\else\global\hoffset=8pt\fi  
  \shipout\vbox{%
    \ifplain@
      \let\makeheadline\relax \let\makefootline\relax
    \else
      \iffirstpage@ \global\firstpage@false
        \let\rightheadline\frheadline
        \let\leftheadline\flheadline
      \else
        \ifrunheads@ 
        \else \let\makeheadline\relax
        \fi
      \fi
    \fi
    \makeheadline \pagebody \makefootline}%
  \advancepageno \ifnum\outputpenalty>-\@MM\else\dosupereject\fi
}
\def\Beta{\mathchar"0\hexnumber@\rmfam 42}
\redefine\mm@{2010} 
\catcode`\@=\active
\nopagenumbers
\chardef\textvolna='176

\chardef\bigalpha='013
\def\negskp{\hskip -2pt}

\chardef\degree="5E
\def\compos{\,\raise 1pt\hbox{$\sssize\circ$} \,}

\def\Span{\operatorname{Span}}

\def\blue#1{#1}

\catcode`#=11\def\diez{#}\catcode`#=6
\catcode`&=11\catcode`&=4
\catcode`_=11\catcode`_=8
\catcode`\^=11\catcode`\^=7
\catcode`~=11\def\volna{~}\catcode`~=\active
\def\mycite#1{\cite{\blue{#1}}\immediate\special{ps:
     ShrHPSdict begin /ShrBORDERthickness 0 def}}
\def\myciterange#1#2#3#4{\cite{\blue{#2#3#4}}\immediate\special{ps:
     ShrHPSdict begin /ShrBORDERthickness 0 def}}
\def\mytag#1{%
    \tag#1}
\def\mythetag#1{\thetag{\blue{#1}}\immediate\special{ps:
     ShrHPSdict begin /ShrBORDERthickness 0 def}}
\def\myrefno#1{\no#1}
\def\myhref#1#2{\blue{#2}\immediate\special{ps:
     ShrHPSdict begin /ShrBORDERthickness 0 def}}
\def\myEarXivlink{\myhref{http://arXiv.org}{http:/\negskp/arXiv.org}}

\def\mytheorem#1{\csname proclaim\endcsname{Theorem #1}}
\def\mytheoremwithtitle#1#2{\csname proclaim\endcsname{Theorem #1#2}}
\def\mythetheorem#1{\blue{#1}\immediate\special{ps:
     ShrHPSdict begin /ShrBORDERthickness 0 def}}
\def\mylemma#1{\csname proclaim\endcsname{Lemma #1}}
\def\mylemmawithtitle#1#2{\csname proclaim\endcsname{Lemma #1#2}}

\def\mycorollary#1{\csname proclaim\endcsname{Corollary #1}}

\def\mydefinition#1{\definition{Definition #1}}
\def\mythedefinition#1{\blue{#1}\immediate\special{ps:
     ShrHPSdict begin /ShrBORDERthickness 0 def}}
\def\myconjecture#1{\csname proclaim\endcsname{Conjecture #1}}
\def\myconjecturewithtitle#1#2{\csname proclaim\endcsname{Conjecture #1#2}}

\def\myproblem#1{\csname proclaim\endcsname{Problem #1}}
\def\myproblemwithtitle#1#2{\csname proclaim\endcsname{Problem #1#2}}


\pagewidth{360pt}
\pageheight{606pt}
\topmatter
\title
On simultaneous approximation of several eigenvalues of a semi-definite 
self-adjoint linear operator in a Hilbert space
\endtitle
\rightheadtext{On simultaneous approximation of several eigenvalues \dots}
\author
Ruslan Sharipov
\endauthor
\address Bashkir State University, 32 Zaki Validi street, 450074 Ufa, Russia
\endaddress
\email
\myhref{mailto:r-sharipov\@mail.ru}{r-sharipov\@mail.ru}
\endemail
\abstract
A lower semi-definite self-adjoint linear operator
in a Hilbert space is taken whose discrete spectrum is not empty 
and comprises at least several eigenvalues  $\lambda_{\text{min}}=\lambda_1\leqslant\ldots\leqslant\lambda_m
<\sigma_{\text{ess}}$. The problem of approximation of these eigenvalues 
by eigenvalues of some linear operator in a finite-dimensional space of 
the dimension $s$ is considered and solved. The accuracy of the approximation 
obtained becomes unlimitedly high as $s\to\infty$. 
\endabstract
\subjclassyear{2010}
\subjclass 35P15, 47A07, 47A75, 47B25, 65D15\endsubjclass
\keywords Hilbert space, self-adjoint operator, semi-definite 
operator, discrete spectrum, eigenvalue, approximation
\endkeywords
\endtopmatter
\loadbold
\TagsOnRight
\document

\head
1. Introduction.
\endhead
     The interest to Hilbert spaces and to linear operators in them is mainly due
to their applications in quantum mechanics where the spaces of square integrable 
complex functions $L_2(\Bbb C,\Bbb R^{3N})$ and differential operators in 
them are considered. The most important of them is the energy operator which is 
called the Hamiltonian in general case and is called the Schr\"odinger operator
in the case of particles in potential fields. In many cases the Schr\"odinger 
operator appears to be self-adjoint and semi-definite. With the quantum mechanical
applications in view, Reed and Simon in their book \mycite{1} give the following 
theorem (see \S~2 in Chapter~\uppercase\expandafter{\romannumeral 13} of \mycite{1}). 
\mytheorem{1.1} Let $F$ be a lower semi-definite self-adjoint linear operator with the 
domain $D(F)$ in a Hilbert space $H$ such that its discrete spectrum is not empty
and $\lambda_{\text{min}}$ is ts minimal eigenvalue\footnotemark.
Let $X$ be an eigenvector of $F$ associated with the eigenvalue $\lambda_{\text{min}}$.
Assume that the eigenvector $X$ is expanded in some orthonormal basis $\{h_i\}_{i=1,\,
\ldots,\,\infty}$ of the Hilbert space $H$:
\footnotetext{\ It is implicitly assumed that $\lambda_{\text{min}}$ is below the
lower limit of the essential spectrum $\sigma_{\text{ess}}(F)$ of the operator $F$.}
$$
\hskip -2em
X=\lim_{n\to\infty}X_n\text{, \ where \ }X_n=\sum^n_{i=1}x^i\,h_i.
\mytag{1.1}
$$
Assume that $h_i\in D(F)$ for all $i=1,\,\ldots,\,\infty$ and suppose that  
$$
\hskip -2em
\exists\ \lim_{n\to\infty}\langle X_n|FX_n\rangle=\lambda_{\text{min}}
\,\Vert X\Vert^2.
$$
Under these assumptions we have the equality  
$$
\hskip -2em
\lambda_{\text{min}}=\lim_{n\to\infty}\hat\mu^{(n)}_1,
$$
where $\hat\mu^{(n)}_1$ is the minimal eigenvalue of the Hermitian $n\times n$ 
matrix $\Phi^{(n)}$ whose elements are $\Phi_{ij}=\langle h_i|Fh_j\rangle$.
\endproclaim
     The angular brackets $\langle\bullet|\bullet\rangle$ in Theorem~\mythetheorem{1.1} 
mean the standard sesquilinear scalar product of the Hilbert space $H$. The formula 
\mythetag{1.1} in Theorem~\mythetheorem{1.1} means convergence of the series 
$$
\hskip -2em
X=\sum^\infty_{i=1}x^i\,h_i,
\mytag{1.2}
$$
with respect to the norm in $H$, while the coefficients $x^i$ in
\mythetag{1.1} and \mythetag{1.2} are denoted according to the Einstein's 
tensorial notation (see \S\,20 of Chapter~\uppercase\expandafter{\romannumeral 1} 
in \mycite{2}), i.\,e. using the upper index $i$ for the coordinates of the vector
$X$.\par
    Theorem~\mythetheorem{1.1} is in background of amost all numerical methods
of quantum chemistry, though it is used without recognition and without references 
to it. The goal of this paper is 
\roster
\rosteritemwd=5pt
\item"1)" to extend Theorem~\mythetheorem{1.1} to the case of several consequtive 
eigenvalues of the operator $F$ taken in the non-decreasing order starting from the lowest 
one;
\item"2)" exclude the orthonormal basis from it;
\item"3)"  replace the domain $D(F)$ of the operator $F$ by the domain $Q(F)$ of 
the corresponding sesquilinear form $q_F(X,Y)=\langle X|FY\rangle$.
\endroster
\par
\head
2. Standard definitions and prerequisites.
\endhead
     A Hilbert space $H$ is a complex (generally speaking infinite-dimensional)
linear vector space with some fixed positive sesquilinear form $\langle X|Y\rangle
=\overline{\langle Y|X\rangle}$ that defines the norm $\Vert X\Vert=
\sqrt{\langle X|X\rangle}$ and thus defines the topology of a complete metric space 
in $H$. With the quantum mechanical applications in view, all Hilbert spaces in this paper
are implicitly assumed to be separable.\par 
     The form $\langle X|Y\rangle$ in $H$ is assumed to be linear in its second argument 
$Y$ and to be conjugate linear with respect to the first argument $X$. This convention 
is used in quantum mechanics (see. \myciterange{3}{3}{, }{4}). The form $\langle X|Y\rangle$ is called 
the standard scalar product in $H$.\par
     A linear operator $F$ in a Hilbert space $H$ is usually given along with its 
domain $D(F)$ which is assumed to be a dense subspace of $H$. The graph of a linear 
operator $F$ is the following subset of the Cartesian product $H\times H$:
$$
\Gamma(F)=\{(X,Y)\in H\times H\!:\ X\in D(F)\text{\ \ and \ }Y=FX\}.
$$
\mydefinition{2.1} 
A linear operator $F$ in a Hilbert space $H$ is called closed, if its 
graph $\Gamma(F)$ is closed in $H\times H$.
\enddefinition
     An extension of an operator $F$ is another operator with the domain bigger than 
$D(F)$ which upon restricting to $D(F)$ coincides with $F$.
\mydefinition{2.2} A linear operator $F$ in a Hilbert space $H$ is called closable
if it has at least one closed extension in $H$. \pagebreak The minimal closed extension 
of a closable operator $F$ is its closed extension with the minimal domain. It is 
denoted through $\bar F$ and is shortly called the closure of $F$. 
\enddefinition
     The minimal closed extension $\bar F$ of a closable operator $F$ is constructed
by closing its graph: $\Gamma(\bar F)=\overline{\Gamma(F)}$. This yields the folloiwing 
theorem.
\mytheorem{2.1} For any closable operator $F$ in a Hilbert space $H$ a vector $X$ 
belongs to the domain $D(\bar F)$ if and only if there is a sequence of vectors 
$X_n\in D(F)$ such that the following relationships hold:
$$
\xalignat 2
&\hskip -2em
\lim_{n\to\infty}X_n=X,
&&\exists\lim_{n\to\infty}FX_n=Y\in H.
\endxalignat
$$
Limits here are in the sense of convergence with respect to the norm in $H$, while
the vector $Y$ is the value of the operator $\bar F$ applied to the vector $X$, i.\,e. 
$\bar FX=Y$.
\endproclaim
\mydefinition{2.3} Let $F$ be a linear operator in a Hilbert space $H$. The linear 
operator $F^*$ with the domain 
$$
\hskip -2em
D(F^*)=\{X\in H\!:\,\exists\,Y\in H\!:\,
\langle X|FZ\rangle=\langle Y|Z\rangle\ \forall\,Z\in D(F)\}
\mytag{2.1}
$$
given by the formula $F^*X=Y$, where the element $Y\in H$ is uniquely fixed by the 
condition $\langle X|FZ\rangle=\langle Y|Z\rangle\ \forall\,Z\in D(F)$ from 
\mythetag{2.1}, is called the conjugate operator for the operator $F$. 
\enddefinition
In the book \mycite{5} the following theorem is given (see \S\,1 of
Chapter~\uppercase\expandafter{\romannumeral8} in \mycite{5}).
\mytheorem{2.2} Let $F$ be a linear operator in a Hilbert space $H$. Then 
\roster
\rosteritemwd=5pt
\item"1)" the conjugate operator $F^*$ is closed;
\item"2)" the operator $F$ is closable if and only if $D(F^*)$ is dense in $H$ and 
$\bar F=F^{**}$ in such a case; 
\item"3)" if the operator $F$ is closable, then $(\bar F)^*=F^*$. 
\endroster
\endproclaim
\mydefinition{2.4} A linear operator $F$ in a Hilbert space $H$ is called symmetric if
the conjugate operator $F^*$ is an extension of $F$.
\enddefinition
Practically the definition~\mythedefinition{2.4} can be replaced with a more simple and
equivalent definition. 
\mydefinition{2.5} A linear operator $F$ in a Hilbert space $H$ is called symmetric
if $\langle X|FY\rangle=\langle FX|Y\rangle$ for all $X,Y\in D(F)$.
\enddefinition
\mytheorem{2.3} Each symmetric operator $F$ in a Hilbert space $H$ is closable. Its 
closure $\bar F$ itself is a symmetric operator in $H$.
\endproclaim
\mydefinition{2.6} A linear operator $F$ in a Hilbert space $H$ is called self-adjoint 
if $F^*=F$, i\.\,e\. if $F$ is symmetric and $D(F)=D(F^*)$.
\enddefinition
\mydefinition{2.7} A symmetric operator $F$ in a Hilbert space $H$ is called
essentially self-adjoint if its closure $\bar F$ is a self-adjoint opertator 
in $H$.
\enddefinition
     Each linear operator $F$ in a Hilbert space $H$ is associated with the sesquilinear 
form $q_F(X,Y)=\langle X|FY\rangle$ which is defined for all $X$ and $Y$ in $D(F)$. In the
case of a symmetric operator $F$ the form $q_F$ is symmetric, i\.\,e\. the following
relationship holds: $q_F(X,Y)=\overline{q_F(Y,X)}\text{\ \ for all \ }X,Y\in D(F)$.
\mydefinition{2.8} A symmetric linear operator $F$ in a Hilbert space $H$ is called 
lower semi-definite if its form $q_F$ is lower semi-definite, i\.\,e\. if there is 
a real constant $C$ such that $q_F(X,X)=\langle X|FX\rangle\geqslant C\,\Vert X\Vert^2
\text{\ \ for all  \ }X\in D(F)$. If $C=0$, such an operator $F$ is called
non-negative. 
\enddefinition
     In the case of a self-adjoint operator $F$, using the spectral theorem, the form
$q_F$ is extended from $D(F)$ to the bigger set 
$$
Q(F)=D(\sqrt{|F|})
\mytag{2.2}
$$
(see \S~6 of Chapter~\uppercase\expandafter{\romannumeral 8} in \mycite{5} and
\S\,53 of Part~\uppercase\expandafter{\romannumeral 15} in \mycite{6}), preserving
its symmetry:
$$
q_F(X,Y)=\overline{q_F(Y,X)}\text{\ \ for all \ }X,Y\in Q(F).
$$
In the case of a lower semi-definite self-adjoint linear operator $F$ the 
set \mythetag{2.2} and the values of the form $q_F$ in it admit a more constructive
description.\par
\mytheorem{2.4} Let $F$ be a lower semi-definite self-adjoint linear operator 
in a Hilbert space $H$. Then a vector $X$ belongs to the domain $Q(F)$ of the 
corresponding sesquilinear form $q_F$ if and only if there is a sequence of vectors 
$X_n$ in the domain $D(F)$ of the operator $F$ such that 
\roster
\rosteritemwd=5pt
\item"1)"\quad $\Vert X_n-X\Vert\to 0$ \ as \ $n\to\infty$;
\item"2)"\quad $\langle X_n-X_m|F(X_n-X_m)\rangle\to 0$
\ as \ $n,m\to\infty$.
\endroster
\endproclaim
\mytheorem{2.5} Under the assumptions of Theorem~\mythetheorem{2.4} let $X$ and $Y$
be two vectors from $Q(F)$ and let $X_n$ and $Y_n$ be two sequences of vectors 
from $D(F)$ approximating $X$ and $Y$ in the sense of Theorem~\mythetheorem{2.4}.
Then
$$
q_F(X,Y)=\lim_{n\to\infty}\langle X_n|FY_n\rangle.
$$
\endproclaim
\head
3. Minimax principle.
\endhead
\mytheorem{3.1} Let $F$ be a lower semi-definite self-adjoint operator with the 
domain $D(F)$ in a Hilbert space $H$. For some positive integer $n$ denote 
$$
\hskip -2em
\mu_n(F)=\sup\limits_{X_1,\,\ldots,\,X_{n-1}}
\ \inf\Sb X\in D(F),\ \Vert X\Vert=1\\
X\perp X_1,\,\ldots,\,X\perp X_{n-1}\endSb\ \langle X|FX\rangle.
\mytag{3.1}
$$
Then exactly one of the two options holds:
\roster
\rosteritemwd=0pt
\item"1)" there exist $n$ eigenvalues of the operator $F$ below the lower limit of its
essential spectrum $\sigma_{\text{ess}}(F)$ and the number $\mu_n(F)$ in \mythetag{3.1}
coincides with $\lambda_n$, provided that each eigenvalue of the operator $F$ in the non-decreasing sequence $$
\lambda_{\text{min}}=\lambda_1\leqslant\ldots\leqslant\lambda_n
$$
is repeated according to its multiplicity, except for maybe the last one;
\item"2)" the number $\mu_n(F)$ in \mythetag{3.1} coincides with the lower limit of 
of the essential spectrum $\sigma_{\text{ess}}(F)$ of the operator $F$. In this case 
there are not more than $n-1$ eigenvalues of the operator $F$ below 
$\sigma_{\text{ess}}(F)$, each being counted with its multiplicity, and 
$\mu_m(F)=\mu_n(F)$ for all $m>n$.
\endroster
\endproclaim
Theorem~\mythetheorem{3.1} expresses the minimax principle for linear operators 
in Hilbert spaces. Its statement and its proof are given in the book \mycite{1} 
(see \S~1 of Chapter~\uppercase\expandafter{\romannumeral13}).\par
     {\bf Remark}. In \mycite{1} it is said that in Theorem~\mythetheorem{3.1}
the domain $D(F)$ of the operator $F$ can be replaced with the domain $Q(F)$ of its
sesquilinear form $q_F$ (see \mythetag{2.2}). At the same time the quantity 
$\langle X|FX\rangle$ is replaced by the value of the form $q_F(X,X)$.
\par
     The case $n=1$ in Theorem~\mythetheorem{3.1} is somewhat different from others. In 
this case the formula \mythetag{3.1} simplifies and takes the form 
$$
\hskip -2em
\mu_1(F)=\inf\limits_{X\in D(F),\ \Vert X\Vert=1}
\langle X|FX\rangle.
\mytag{3.2}
$$
The next theorem is a corollary of Theorem~\mythetheorem{3.1}.
\mytheorem{3.2} If a lower semi-definite self-adjoint operator $F$ with the domain 
$D(F)$ in a Hilbert space $H$ has a non-empty discrete spectrum below 
$\sigma_{\text{ess}}(F)$, then its minimal eigenvalue $\lambda_{\text{min}}=\mu_1(F)$ 
is given by the formula \mythetag{3.2}.
\endproclaim
{\bf Remark}. The domain $D(F)$ of the operator $F$ in the formula \mythetag{3.2}
can be replaced by the domain $Q(F)$ of its sesquilinear form $q_F$ (see \mythetag{2.2}). 
At the same time the quantity $\langle X|FX\rangle$ is replaced by the value of the
form $q_F(X,X)$.\par 
\head
4. Minimum principle.
\endhead
     In the case $n=1$ the minimax principle is expressed by the formula
\mythetag{3.2}. It transforms into minimum principle. A similar minimum principle can 
be formulated for the case $n>1$. 
\mytheorem{4.1} Let $F$ be a lower semi-definite self-adjoint linear operator with the
domain $D(F)$ in a Hilbert space $H$ whose discrete spectrum below $\sigma_{\text{ess}}(F)$
comprises at least $n$ eigenvalues arranged in the non-decreasing order
$$
\lambda_{\text{min}}=\lambda_1\leqslant\ldots\leqslant\lambda_{n-1}\leqslant\lambda_n
$$ 
so that each eigenvalue in the sequence is repeated according to its multiplicity, 
except for maybe the last one. Let $X_1,\,\ldots,\,X_{n-1}$ be linerly independent 
eigenvectors for the initial $n-1$ eigenvalues in this non-decreasing sequence. 
Then the quantity $\mu_n(F)=\lambda_n$ is given by the formula 
$$
\mu_n(F)=\inf\Sb X\in D(F),\ \Vert X\Vert=1\\
X\perp X_1,\,\ldots,\,X\perp X_{n-1}\endSb
\langle X|FX\rangle.
$$
\endproclaim
     The minimum principle similar to Theorem~\mythetheorem{4.1} for symmetric operators 
in a finite-dimensional Euclidean space can be found in the book \mycite{7}. For the Laplace
operator it is formulated in \mycite{8}. The proof of Theorem~\mythetheorem{4.1} can be found
in \mycite{9}. In this paper Theorem~\mythetheorem{4.1} is given for the sake of completeness 
of our preliminary review. It is not used in what follows.\par
\head
5. The Rayleigh-Ritz method.
\endhead
\mytheorem{5.1} Let $F$ be a lower semi-definite self-adjoint linear operator with the
domain $D(F)$ in a Hilbert space $H$, let $V\subset D(F)$ be a finite-dimensional
subspace of the dimension $n$ in $D(F)$, and let $P$ be the orthogonal projector onto
$V$. The composite operator $P\compos F\compos P$ is lower semi-definite and self-adjoint.
Its restriction $F_V$ to the subspace $V$ has exactly $n$ eigenvalues $\hat\mu_1,\,\ldots,\,\hat\mu_n$ that can be arranged in the non-decreasing order 
$\hat\mu_1\leqslant\ldots\leqslant\hat\mu_n$ where each eigenvalue is repeated according
to its multiplicity. Under these assumptions the following inequalities hold:
$$
\hskip -2em
\mu_m(F)\leqslant\hat\mu_m\text{, \ where \ }m=1,\,\ldots,\,n.
\mytag{5.1}
$$
\endproclaim
     The quantities $\mu_m(F)$ in the left hand side of the inequalities 
\mythetag{5.1} are given by the formula \mythetag{3.1}. They do not depend on the 
choice of the subspace $V$ in Theorem~\mythetheorem{5.1}. The statement and the 
proof of Theorem~\mythetheorem{5.1} are given in the book \mycite{1} (see \S~2 of 
Chapter~\uppercase\expandafter{\romannumeral 13}).\par
     {\bf Remark}. As we already noted above, the domain $D(F)$ of the operator $F$ 
when computing $\mu_m(F)$ in the formula \mythetag{3.1} can be replaced by the
domain $Q(F)$ of its sesquilinear form $q_F$ (see \mythetag{2.2}). At the same time 
the quantity $\langle X|FX\rangle$ is replaced by the value of the form $q_F(X,X)$.
In Theorem~\mythetheorem{5.1} the domain $D(F)$ can also be replaced by the domain 
$Q(F)$. Although then we should define the operator $F_V$ in a different way 
through the restriction of the form $q_F$ to the subspace $V\subset Q(F)$. 
\mytheorem{5.2} Let $F$ be a lower semi-definite self-adjoint linear operator
in a Hilbert space $H$ whose associated sesquilinear form is $q_F$ and the domain 
of $q_F$ is $Q(F)$. Let $V\subset Q(F)$ be a finite-dimensional
subspace of the dimension $n$ in $Q(F)$ and let $F_V$ be a linear operator in the 
subspace $V$ defined by the restriction of the form $q_F$ to $V$ according to
the formula
$$
\langle Y|F_VX\rangle=q_F(Y,X)\text{\ \ for all \ }X,Y\in V.
$$
The operator $F_V$ has exactly $n$ eigenvalues $\hat\mu_1,\,\ldots,\,\hat\mu_n$ that 
can be arranged in the non-decreasing order $\hat\mu_1\leqslant\ldots\leqslant\hat\mu_n$ 
where each eigenvalue is repeated according to its multiplicity. Under these assumptions 
the following inequalities hold:
$$
\hskip -2em
\mu_m(F)\leqslant\hat\mu_m\text{, \ where \ }m=1,\,\ldots,\,n.
\mytag{5.2}
$$
\endproclaim
\demo{Proof} Applying the minimax principle expressed by Theorem~\mythetheorem{3.1}
to the operator $F_V$ in the finite-dimensional space $V$, from \mythetag{3.1} we 
derive
$$
\hskip -2em
\gathered
\hat\mu_m=\sup\limits_{X_1,\,\ldots,\,X_{m-1}\in V}
\ \inf\Sb X\in V,\ \Vert X\Vert=1\\
X\perp X_1,\,\ldots,\,X\perp X_{m-1}\endSb\ \langle X|F_VX\rangle=\\
=\sup\limits_{X_1,\,\ldots,\,X_{m-1}\in V}
\ \inf\Sb X\in V,\ \Vert X\Vert=1\\
X\perp X_1,\,\ldots,\,X\perp X_{m-1}\endSb\ q_F(X,X).
\endgathered
\mytag{5.3}
$$
Let $P$ be the orthogonal projector onto the subspace $V$. Them \mythetag{5.3}
implies  
$$
\hskip -2em
\gathered
\hat\mu_m=\sup\limits_{X_1,\,\ldots,\,X_{m-1}\in H}
\ \inf\Sb X\in V,\ \Vert X\Vert=1\\
X\perp PX_1,\,\ldots,\,X\perp PX_{m-1}\endSb\ q_F(X,X)=\\
=\sup\limits_{X_1,\,\ldots,\,X_{m-1}\in H}
\ \inf\Sb X\in V,\ \Vert X\Vert=1\\
X\perp X_1,\,\ldots,\,X\perp X_{m-1}\endSb\ q_F(X,X)\geqslant\\
\geqslant\sup\limits_{X_1,\,\ldots,\,X_{m-1}\in H}
\ \inf\Sb X\in Q(F),\ \Vert X\Vert=1\\
X\perp X_1,\,\ldots,\,X\perp X_{m-1}\endSb\ q_F(X,X)=\mu_m(F).
\endgathered
\mytag{5.4}
$$
The inequality in \mythetag{5.4} arises since we replace $V$ in the left hand 
side of the inequality by a larger subspace $Q(F)$ in the right hand side of 
this inequality. The last equality in \mythetag{5.4} is due to the remark to Theorem~\mythetheorem{3.1} on page 5.
\qed\enddemo
     Theorems~\mythetheorem{5.1} and \mythetheorem{1.1} constitute a base for 
the Rayleigh-Ritz method. Its application to the experimental confirmation of 
the Lamb shift is described in \mycite{1} (see \S~2 
of Chapter~\uppercase\expandafter{\romannumeral 13} therein). See also
\myciterange{10}{10}{--}{17}.\par 
\head
6. Approximation of several eigenvalues.
\endhead
\mytheorem{6.1} Let $F$ be a lower semi-definite self-adjoint linear operator with 
the domain $D(F)$ in a Hilbert space $H$ whose discrete spectrum below 
$\sigma_{\text{ess}}(F)$ is not empty and comprises at least $m$ eigenvalues 
arranged in the non-decreasing order $\lambda_{\text{min}}=\lambda_1\leqslant\ldots
\leqslant\lambda_m$ so that each eigenvalue is repeated according to its multiplicity,
except for maybe the last one. Let $X_1,\,\ldots,\,X_m$ be linearly independent
eigenvectors corresponding to the eigenvalues $\lambda_{\text{min}}=\lambda_1\leqslant
\ldots\leqslant\lambda_m$ and assume that for each eigenvector $X_k$ a sequence of 
vectors $X_{kn}$ from the domain $Q(F)$ of the associated with $F$ sesquilinear form 
$q_F$ is given so that
$$
\hskip -2em
X_k=\lim_{n\to\infty} X_{kn}\text{\ \ and \ }
\lim_{n\to\infty}q_F(X_{kn},X_{qn}\rangle=
\langle X_k|FX_q\rangle,
\mytag{6.1}
$$
where $1\leqslant k,q\leqslant m$. Under these assumptions 
$$
\hskip -2em
\lambda_k=\lim_{n\to\infty}\hat\mu^{(n)}_k\text{\ \ for all \ }k=1,\,\ldots,m,
\mytag{6.2}
$$
where $\hat\mu^{(n)}_{\text{min}}=\hat\mu^{(n)}_1\leqslant\ldots\leqslant
\hat\mu^{(n)}_m$ are initial $m$ eigenvalues of the operator $F_{V_n}$ generated by 
the restriction of the form $q_F$ to the finite-dimensional subspace
$$
\hskip -2em
V_n=\Span(\{X_{ks}\text{, \ where \ }1\leqslant k\leqslant m,
1\leqslant s\leqslant n\})\subset Q(F)
\mytag{6.3}
$$
taken in the non-decreasing order so that each eigenvalue is repeated according to its multiplicity, except for maybe the last one. 
\endproclaim
\demo{Proof} The eigenvectors $X_1,\,\ldots,\,X_m$ belong to the domain $D(F)$. They are
linearly independent. Their span 
$$
\hskip -2em
\Cal V=\Span(X_1,\ldots,\,X_m)\subset D(F)\subset Q(F)
\mytag{6.4}
$$
is an $m$-dimensional subspace in $D(F)$ and in $Q(F)$. Let $\Cal P$ be the orthogonal
projector onto the subspace \mythetag{6.4}. We set 
$$
\hskip -2em
F_{\Cal V}=\Cal P\compos F\compos \Cal P\ \hbox{\vrule 
height 10pt depth 8pt width 0.5pt}_{\,\,\Cal V}.
\mytag{6.5}
$$
In other words, we denote through $F_{\Cal V}$ the restriction of the composite operator 
$\Cal P\compos F\compos\Cal P$ to the subspace \mythetag{6.4}. The eigenvectors 
$X_1,\,\ldots,\,X_m$ of the operator $F$ are eigenvectors of the operator 
\mythetag{6.5}, while the corresponding eigenvalues coincide with $\lambda_{\text{min}}
=\lambda_1\leqslant\ldots\leqslant\lambda_m$. Therefore the matrix of the operator 
\mythetag{6.5} in the basis of the eigenvectors $X_1,\,\ldots,\,X_m$ is 
$$
\hskip -2em
\Cal F=\Vmatrix 
\lambda_1 &\hdots &0\\
\vdots & \ddots & \vdots\\
0 &\hdots &\lambda_m
\endVmatrix.
\mytag{6.6}
$$
The operator \mythetag{6.5} is generated by the restriction of the form $q_F$ to the
subspace \mythetag{6.4}, therefore the matrix \mythetag{6.6} obeys the relationships 
$$
\xalignat 2
&\hskip -2em
\Cal H=G\,\Cal F,
&&\Cal F=G^{-1}\,\Cal H,
\mytag{6.7}
\endxalignat
$$
where $\Cal H$ and $G$ are two Hermitian matrices with the components 
$$
\xalignat 2
&\hskip -2em
\Cal H_{sk}=q_F(X_s,X_k)=\langle X_s|FX_k\rangle,
&&g_{sq}=\langle X_s|X_q\rangle.
\mytag{6.8}
\endxalignat
$$
The matrix $G$ is the Gram matrix of the basis of vectors $X_1,\,\ldots,\,X_m$. 
Since these vectors are linearly independent, $\det G\neq 0$ and hence the matrix 
relationships \mythetag{6.7} are consistent.\par
     Let's consider the vectors $X_{1n},\,\ldots,\,X_{mn}$. According to the first
condition in \mythetag{6.1} these vectors approximate $X_1,\,\ldots,\,X_m$ in the sense
of convergence with respect to the norm in the Hilbert space $H$. This means that
for all sufficiently large $n$ the vectors $X_{1n},\,\ldots,\,X_{mn}$ are linearly
independent. We set 
$$
\hskip -2em
\Cal V_n=\Span(X_{1n},\ldots,\,X_{mn})\subset Q(F).
\mytag{6.9}
$$
The subspace \mythetag{6.9} is analogous to the subspace \mythetag{6.4}. For all 
sufficiently large $n$ the dimensions of these subspaces do coincide:  
$$
\hskip -2em
\dim\Cal V_n=\dim\Cal V=m.
\mytag{6.10}
$$
Generally speaking, the subspace \mythetag{6.9} is not enclosed in $D(F)$. Therefore
we cannot write a formula similar to \mythetag{6.5}. But we can consider the operator
$F_{\Cal V_n}$ generated by the restriction of the sesquilinear form $q_F$ to the
subspace \mythetag{6.9}. Let $\Cal F^{(n)}$ be the matrix of such operator in the basis 
$X_{1n},\,\ldots,\,X_{mn}$. It is similar to the matrix \mythetag{6.6}, though it is not
diagonal. The matrix $\Cal F^{(n)}$ obeys the relationships which are similar to
the relationships \mythetag{6.7}: 
$$
\xalignat 2
&\hskip -2em
\Cal H^{(n)}=G^{(n)}\,\Cal F^{(n)},
&&\Cal F^{(n)}=(G^{(n)})^{-1}\,\Cal H^{(n)}.
\mytag{6.11}
\endxalignat
$$
In \mythetag{6.11} we see two Hermitian matrices $\Cal H^{(n)}$ and $G^{(n)}$. Their
components are defined by the formulas similar to \mythetag{6.8}:
$$
\xalignat 2
&\hskip -2em
\Cal H[n]_{sk}=q_F(X_{sn},X_{kn})=\langle X_{sn}|F_{\Cal V_n}X_{kn}\rangle,
&&g[n]_{sq}=\langle X_{sn}|X_{qn}\rangle.
\mytag{6.12}
\endxalignat
$$
The matrix $G^{(n)}$ in \mythetag{6.11} is the Gram matrix for the basis composed 
by the vectors $X_{1n},\,\ldots,\,X_{mn}$. For all sufficiently large $n$ the matrix
$G^{(n)}$ is non-degenerate, which is in agreement with \mythetag{6.10}. Hence the relationships \mythetag{6.11} are consistent.\par
     Now we consider the conditions \mythetag{6.1}. Applying them to \mythetag{6.8} 
and \mythetag{6.12}, we get 
$$
\xalignat 2
&\hskip -2em
\lim_{n\to\infty}g[n]_{sq}=g_{sq}, 
&&\lim_{n\to\infty}\Cal H[n]_{sk}=\Cal H_{sk}. 
\mytag{6.13}
\endxalignat
$$
In the matrix form the relationships \mythetag{6.13} are written as follows:
$$
\xalignat 2
&\hskip -2em
\lim_{n\to\infty}G^{(n)}=G, 
&&\lim_{n\to\infty}\Cal H^{(n)}=\Cal H. 
\mytag{6.14}
\endxalignat
$$
The first relationship \mythetag{6.14} yields $\lim\limits_{n\to\infty}\det G^{(n)}
=\det G\neq 0$. Hence
$$
\hskip -2em
\exists\ \lim_{n\to\infty}(G^{(n)})^{-1}=G^{-1}. 
\mytag{6.15}
$$
Let's combine \mythetag{6.15} with the second relationship \mythetag{6.14} and then 
take into account \mythetag{6.11}. This yields the following relationship
$$ 
\hskip -2em
\exists\ \lim_{n\to\infty}\Cal F^{(n)}=\Cal F. 
\mytag{6.16}
$$
The eigenvalues of the matrix $\Cal F$ are the eigenvalues of the operator \mythetag{6.5}, coinciding with $\lambda_{\text{min}}=\lambda_1\leqslant\ldots\leqslant\lambda_m$. They are
presented explicitly in \mythetag{6.6}. The eigenvalues of the matrix $\Cal F^{(n)}$ in \mythetag{6.16} are the eigenvalues of the self-adjoint operator $F_{\Cal V_n}$ in 
$m$-dimensional space $\Cal V_n$. We denote them 
$$
\hskip -2em
\hat\lambda_{\text{min}}^{(n)}=\hat\lambda_1^{(n)}\leqslant\ldots
\leqslant\hat\lambda_m^{(n)}.
\mytag{6.17}
$$
From \mythetag{6.16} for the eigenvalues \mythetag{6.17} we derive 
$$
\hskip -2em
\exists\ \lim_{n\to\infty}\hat\lambda_k^{(n)}=\lambda_k
\text{\ \ for all \ }k=1,\,\ldots,m. 
\mytag{6.18}
$$
\par
     In the next step we consider the subspace \mythetag{6.3}. The sesquilinear form 
$q_F$ produces the self-adjoint operator $F_{V_n}$ whose initial $m$ eigenvalues
in Theorem~\mythetheorem{6.1} are denoted through $\hat\mu^{(n)}_{\text{min}}=\hat
\mu^{(n)}_1\leqslant\ldots\leqslant\hat\mu^{(n)}_m$. We can apply 
Theorem~\mythetheorem{5.2} to the subspace \mythetag{6.3}. Applying this theorem, 
we get 
$$
\hskip -2em
\mu_k(F)\leqslant\hat\mu_k\text{, \ where \ }k=1,\,\ldots,\,m.
\mytag{6.19}
$$
Moreover, in our case the operator $F$ is such that for all $\mu_k(F)$ from 
\mythetag{6.19} the first of the two alternative options from Theorem~\mythetheorem{3.1} 
is realized. This means that 
$$
\hskip -2em
\mu_k(F)=\lambda_k\text{ \ for all \ }k=1,\,\ldots,\,m.
\mytag{6.20}
$$
Comparing \mythetag{6.19} and \mythetag{6.20}, we derive
$$
\hskip -2em
\lambda_k\leqslant\hat\mu_k\text{ \ for all \ }k=1,\,\ldots,\,m.
\mytag{6.21}
$$
\par
     In order to complete the proof of Theorem~\mythetheorem{6.1} we compare the
subspaces \mythetag{6.3} and \mythetag{6.9} along with the self-adjoint operators 
$F_{V_n}$ and $F_{\Cal V_n}$ in them. It is easy to see that $\Cal V_n\subset V_n$. 
If we denote through $\Cal P_n$ the orthogonal projector onto the smaller subspace 
$\Cal V_n$, we can easily derive the relationship
$$
\hskip -2em
F_{\Cal V_n}=\Cal P_n\compos F_{V_n}\compos\Cal P_n\ \hbox{\vrule 
height 10pt depth 8pt width 0.5pt}_{\,\,\ \Cal V_n}. 
\mytag{6.22}
$$
Due to the inclusion $\Cal V_n\subset V_n$ and the relationship \mythetag{6.22} 
we can apply Theorem~\mythetheorem{5.1} to the subspace $\Cal V_n$ and to the
operator $F_{V_n}$ in the enclosing subspace $V_n$. Applying this theorem, we
derive the inequalities 
$$
\hskip -2em
\mu_k(F_{V_n})=\hat\mu_k\leqslant\hat\lambda_k^{(n)}
\text{ \ for all \ }k=1,\,\ldots,\,m.
\mytag{6.23}
$$
From \mythetag{6.21} and \mythetag{6.23} we derive double inequalities:
$$
\hskip -2em
\lambda_k\leqslant\hat\mu_k\leqslant\hat\lambda_k^{(n)}
\text{ \ for all \ }k=1,\,\ldots,\,m.
\mytag{6.24}
$$
Now, applying \mythetag{6.18} to \mythetag{6.24}, we get the required result
expressed by the formula \mythetag{6.2}. Theorem~\mythetheorem{6.1} is proved.
\qed\enddemo
\head
7. Conclusions 
\endhead
     The main result of this paper is expressed by Theorem~\mythetheorem{6.1}.
As a mathematical result it has independent value. As for its application to 
quantum mechanics and quantum chemistry, it means that the Rayleigh-Ritz method 
can be used not only for computing the ground energy level of atoms and molecules, 
but for several consequtive excited energy levels along with the ground level. 
A way of applying Theorem~\mythetheorem{6.1} to quantum chemistry using 
polylinear splines is described in \mycite{9}.
\par 

\enddocument
\end